\documentclass[12pt]{article}
\usepackage[latin1]{inputenc}
\usepackage{amssymb,amsmath,amsfonts}
\newtheorem{theorem}{Theorem}
\newtheorem{proposition}[theorem]{Proposition}
\newtheorem{lemma}[theorem]{Lemma}

\newtheorem{corollary}[theorem]{Corollary}

\newcommand{\R}{\mathbb{R}}

\newcommand{\Q}{\mathbb{Q}}
\newcommand{\Sf}{\mathbb{S}}

\newcommand{\Hy}{\mathbb{H}}
\newcommand{\spa}{\mbox{span}}

\newcommand{\grad}{\mbox{grad}}

\newcommand{\po}{{\hspace*{-1ex}}{\bf .  }}
\newcommand{\ii}{isometric immersion }

\newcommand{\Sp}{\mathbb{S}}
\newcommand{\nat}{\nabla^\theta}

\newcommand{\nab}{\tilde\nabla}

\def\Ral{{\cal R}}

\def\<{\langle}

\def\>{\rangle}
\def\a{\alpha}

\def\bea{\begin{eqnarray*} }
\def\eea{\end{eqnarray*} }
\def\be{\begin{equation} }
\def\ee{\end{equation} }

\def\nap{\nabla^\perp}
\def\proof{\noindent{\it Proof: }}
\def\qed{\ifhmode\unskip\nobreak\fi\ifmmode\ifinner
\else\hskip5 pt \fi\fi\hbox{\hskip5 pt \vrule width4 pt
height6 pt  depth1.5 pt \hskip 1pt }}
\begin{document}

\title{The associated family of an  elliptic surface and\\ an application
to minimal submanifolds}
\author{Marcos Dajczer and Theodoros Vlachos}
\date{}
\maketitle

\begin{abstract} It is  well-known  that in any codimension a simply connected 
Euclidean minimal surface has an associated one-parameter family 
of minimal  isometric deformations. In this paper, we show that this is just a special case of the  
associated family to any simply connected elliptic surface for which all curvature 
ellipses of a certain order are circles.  We also provide the conditions under which 
this associated  family is trivial, extending the known result for  minimal surfaces.
As an application, we show how the associated family of a minimal Euclidean submanifold 
of rank two is determined by the associated family of an elliptic surface clarifying 
the geometry around the associated family of these higher dimensional 
submanifolds.
\end{abstract}

\section{Introduction}

It is a well-known fact that a simply connected minimal surface in a space form  
of any dimension allows a one-parameter family of isometric  minimal deformations,
called the associated family, and that in Euclidean space this family can be 
parametrically given by means of the generalized Weierstrass representation; 
see \cite{HO}.  In this paper, we show that this associated family is just a 
special case of the associated family to an elliptic surface for which all
ellipses of curvature of a certain order are circles. Minimal surfaces can 
be seen as those elliptic surfaces for which the ellipse  of curvature of
order zero is  a circle. Several basic properties of the new associated family
are also given, in particular, we state when the family is trivial. 
Our second main result is an application of the result on surfaces to 
minimal submanifolds, which was our initial motivation and is explained in 
the sequel.

 Euclidean  submanifolds of rank two have been studied in different contexts; see 
\cite{BKV}, \cite{DF1}, \cite{DF2}, \cite{DM1} and \cite{DM2}. A submanifold having 
rank two means that the image of the Gauss map is a surface in the corresponding 
Grassmannian or, equivalently, that the kernel of the second fundamental form 
(relative nullity subspace) has constant codimension two. The study of the  minimal ones
is particularly interesting since they belong to the important class of austere 
submanifolds introduced in \cite{HL}. As a special case, one has the ones that carry a
Kaehler structure  described in \cite{DF2} by a  Weierstrass type
representation in terms of $m$-isotropic surfaces. A minimal surface is called 
$m$-isotropic if all ellipses of curvature  up to order $m$ are circles. In turn, the
$m$-isotropic surfaces can be constructed  by the use of a Weierstrass type representation 
given in \cite{DG2} based on results in \cite{Ch}.

It turns that the normal bundle of an elliptic surface splits as the orthogonal sum of 
a sequence of plane bundles (except the last one in odd codimension) such that each fiber 
contains  an ellipse of curvature that is then ordered accordingly; see next section
for details.  Euclidean minimal submanifolds of rank two  have been parametrically 
described in  \cite{DF2} by means of the class of elliptic surfaces for which the 
curvature ellipses of a certain order are  circles.  

To some surprise, it was observed in \cite{DG} that any simply connected Euclidean
minimal submanifold of rank two allows an associated family of submanifolds 
of the same class. As in the surface case, this family 
is obtained  by rotating the second fundamental form while keeping fixed the normal bundle 
and the induced normal connection. This fact, together with the representation in \cite{DF2} 
discussed above, suggests that an elliptic surface in a space form for which the
ellipses of curvature of a certain order are circles should have some kind of associate 
family preserving that property, and that was the starting point  of this paper. 

There is an abundance of examples of surfaces with circular ellipses of curvature, 
specially minimal ones. In particular, there are the surfaces for which all but the 
last one ellipse  of curvature is a circle. These have been studied in the round sphere  
\cite{BPW} and in hyperbolic space \cite{Hu} under the name of superconformal.  
Other interesting examples are holomorphic curves in the nearly Kaehler sphere $\Sf^6$.
The theory of these surfaces started in  \cite{Br} and  was developed in \cite{BVW}, 
\cite{H1} and \cite{H2}. The first ellipse of curvature is always a circle but there 
is a class for which second  curvature ellipse is not a circle; 
see Case $3$ of Theorem 6.5 in \cite{H1}.

For the purpose of this paper, the most important known examples are Lawson's surfaces. 
These are minimal surfaces in spheres  that decompose as a direct sum of elements 
in the associated family 
$h_\theta,\,\theta\in [0,\pi)$, of a minimal surface $h$ in $\Sf^3$. 
More precisely, we consider surfaces in $\Sf^{4n-1}\subset\R^{4n}$ given as
$$
f=a_1h_{\theta_1}\oplus\ldots\oplus a_nh_{\theta_n}
$$
where $0\leq\theta_1<\cdots<\theta_n<\pi$, the  real numbers $a_1,\ldots,a_n$  
satisfy $\sum_{j=1}^na_j^2=1$ and $\oplus$ denotes the orthogonal sum with 
respect to an orthogonal decomposition of $\R^{4n}$.
It has been checked in \cite{Vl} that all ellipses of curvature of even order are
circles while that the ones of odd order generically are not.
These surfaces are part of Lawson's conjecture \cite{La} which
asserts that the only non-flat minimal surfaces in spheres that are locally
isometric to minimal surfaces in $\Sf^3$ are Lawson's surfaces.

Most of what is done in this paper for surfaces can be extended to elliptic submanifolds
of rank two. But in the final section of the paper, we limit ourselves to show how the 
associated family of a minimal Euclidean submanifold of rank two is determined by
the associated family to an elliptic surface with a circular ellipse of
curvature. This result completely clarifies the geometry around the associated family
of these higher dimensional submanifolds.

Finally, we observe that a key ingredient of our proofs is the classical 
Burstin-Mayer-Allendoerfer theory as discussed in Vol.\ IV of Spivak \cite{Sp}. 
Similar to the case of curves, this theory shows that certain tensors associated to 
a set of Frenet type equations are a complete set of invariants for a submanifold of 
a space form. Among these tensors, one has the higher order fundamental forms some of 
which are preserved by our associated family. We should point out that isometric deformations of 
submanifolds that also preserve higher fundamental forms, starting with the second 
fundamental form, up to a stated order was somehow considered in \cite{Bo}.

\section{Preliminaries}

In this section we  recall from \cite{Sp} some basic definitions for submanifolds 
in space forms,  and from \cite{DF2} the notions of elliptic surface, ellipse of 
curvature and polar surface to an elliptic surface and some of their basic properties,
which will be used in the sequel without further reference.
\vspace{1,5ex} 

Let $f\colon M^n\to\Q_{c}^N$ be a substantial \ii of a connected 
$n$-dimensional Riemannian manifold into either the Euclidean space
$\R^N$ ($c=0$), the round sphere $\Sf^N$ ($c>0$) or the hyperbolic 
space $\Hy^N$ ($c<0$) with vector valued second fundamental form $\alpha_f$ 
and induced Riemannian connection$\nap$ in the normal bundle $N_fM$.  
That $f$ is \emph{substantial} (called full in \cite{Sp}) means that the 
codimension cannot be reduced.

The $k^{th}$\emph{-normal space} $ N^f_k(x)$ of $f$ at $x\in M^n$ 
for $k\ge 1$ is defined as
$$
N^f_k(x)=\spa\{\alpha_f^{k+1}(X_1,\ldots,X_{k+1}):X_1,\ldots,X_{k+1}\in T_xM\}.
$$
Thus $\alpha_f^2=\alpha_f$ and for $s\geq 3$  the symmetric tensor
$\alpha_f^s\colon TM\times\cdots\times TM\to N_fM$,  called
the $s^{th}$\emph{-fundamental form}, is defined inductively by
$$
\alpha_f^s(X_1,\ldots,X_s)=\left(\nabla^\perp_{X_s}\ldots
\nabla^\perp_{X_3}\alpha_f(X_2,X_1)\right)^\perp
$$
where $(\;\;)^\perp$ denotes taking the projection onto
the normal subspace $(N^f_1\oplus\ldots\oplus N^f_{s-2})^\perp$.

We always admit that $f$ is \emph{regular} (called nicely curved in \cite{Sp})
which means that all the $N^f_k$'$s$ have constant dimension for each $k$ and thus
form normal subbundles.  This means ``geometrically'' that at each point the submanifold bends 
in the same number of directions. For any submanifold this condition is verified 
along connected components of an open dense subset of $M^n$. 

A surface $g\colon L^2\to\Q_c^N$ is called \textit{elliptic} in \cite{DF2} if there exists a 
(unique up to a sign) almost complex structure $J\colon TL\to TL$ such that 
the second fundamental form satisfies
$$
\alpha_g (X,X)+ \alpha_g (JX,JX)=0 \;\;\mbox{for all}\;\;  X \in TL.
$$
Then all the $N^g_k$'$s$ have dimension two except the last one that is one-dimensional 
if the codimension is odd. Therefore, the normal bundle $N_gL$  splits as 
$$
N_gL=N^g_1\oplus \cdots \oplus N^g_{\tau},
$$
where  $\tau$ is the 
index of the last subbundle. Thus, the induced bundle $g^*T\Q_c^N$ splits as 
$$
g^*T\Q_c^N=N^g_0\oplus N^g_1\oplus \cdots \oplus N^g_{\tau}
$$
where $N^g_0=g_*TL$.  Setting
$$
\tau^o = \left\{\begin{array}{l}
\tau\;\;\;\;\;\;\;\;\;\;\;\mbox{if}\;\;N\;\;\; \mbox{is even}\\
\tau-1\;\;\;\;\; \mbox{if}\;\;N\;\;\; \mbox{is odd},
\end{array} \right.
$$ 
it turns out that the almost complex structure $J$ on $TL$ induces an almost
complex  structure $J_s$ on each $N_s^g$, $1\leq s\leq\tau^o$, defined by
$$
J_s\alpha^{s+1}_g(X_1,\ldots,X_{s+1})=\alpha^{s+1}_g(JX_1,\ldots,X_{s+1}).
$$
In the sequel, we denote by $\pi_s\colon  g^*T\Q_c^N\to N_s^g,\;0\leq s\leq \tau$, the 
orthogonal projection. Then, we have for $2\leq s\leq\tau^o$ that
\be\label{js}
J_s\pi_s(\nap_X\xi)=\pi_s(\nap_XJ_{s-1}\xi)=\pi_s(\nap_{JX}\xi)
\;\;\mbox{if}\;\;\xi\in N_{s-1}^g
\ee
and 
\be\label{jss}
J^t_{s-1}\pi_{s-1}(\nap_X\xi)=\pi_{s-1}(\nap_XJ^t_s\xi)=
\pi_{s-1}(\nap_{JX}\xi)\;\;\mbox{if}\;\;\xi\in N_s^g.
\ee
For any $\varphi\in\Sp^1=[0,\pi)$  let $R^s_\varphi\colon N^g_s\to N^g_s$, 
$0\leq s\leq\tau^o$, denote the map given by  
\be\label{rotation}
R^s_\varphi=\cos\varphi I+\sin\varphi J_s.
\ee
It follows from (\ref{js}) and (\ref{jss}) that
\be\label{one0}
R^{s+1}_\varphi\pi_{s+1}(\nap_X\xi)=\pi_{s+1}(\nap_XR^s_\varphi\xi)
\;\;\mbox{if}\;\;\xi\in N_s^g
\ee
and 
\be\label{two0}
(R^s_\varphi)^t\pi_s(\nap_X\xi)=\pi_s(\nap_X(R^{s+1}_\varphi)^t\xi)
\;\;\mbox{if}\;\;\xi\in N_{s+1}^g
\ee
for any $1\leq s\leq\tau^o-1$.
\medskip

The \textit{$s^{th}$-order curvature ellipse} 
$\mathcal{E}_s^g(x)\subset N^g_s(x)$ of $g$ at $x\in L^2$ for $0\leq s\leq\tau^o$ is 
$$
\mathcal{E}_s^g(x)=\{\alpha_g^{s+1}(Z_{\psi},\dots,Z_{\psi}) : 
Z_{\psi}=\cos\psi Z+\sin\psi JZ\;\;\mbox{and}\;\;\psi\in [0,\pi)\},
$$
where we understand that $\a^1=g_*$ and assume that $Z\in T_xL$ has  unit length 
and satisfies $\<Z,JZ\>=0$. 
It follows from the ellipticity condition that such a $Z$ always exists and that
$\mathcal{E}_s^g(x)$ is indeed  an ellipse. 
\medskip

We point out that $\mathcal{E}_1^g$ given by the above 
definition coincides with the standard definition only if the mean 
curvature vanishes, in which case the higher order ellipses also coincide.
\medskip

By $\mathcal{E}_\ell^g$ \emph{being a circle} we mean that the curvature ellipse 
$\mathcal{E}_\ell^g(x)$ is a circle for any $x\in L^2$. A fundamental fact in 
this paper is that $\mathcal{E}_s^g(x)$ is a circle if and only if $J_s(x)$ is 
orthogonal. Notice that  $\mathcal{E}_0^g$ is a circle if and only if $g$ 
is a minimal surface.
\medskip

A \emph{polar surface} to an elliptic surface 
$g\colon L^2\to\Q_c^{N-c}\subset\R^N$ $(c=0,1)$ is an immersion defined  
as follows:
\begin{itemize}
\item [(i)] If $N-c$ is odd, then the polar surface $h\colon L^2\to\Sf_1^{N-1}$ 
is the spherical image of a unit normal field spanning the last one-dimensional
normal bundle.
\item [(ii)] If $N-c$ is even, then  the polar surface $h\colon L^2\to\R^N$ is 
any surface such that $T_{h(x)}L=N^g_{\tau}(x)$ up to parallel 
identification in $\R^N$.
\end{itemize}

It is known that in case $(ii)$  any elliptic surface admits locally many polar 
surfaces. 
It turns out that a polar surface to an elliptic surface is necessarily elliptic.  
Moreover, if the elliptic surface has a circular ellipse of curvature then its 
polar surface has the same property at the ``corresponding" normal bundle. 
In particular, for the  polar surface to an $m$-isotropic  surface the last 
$m+1$ ellipses  of curvature are circles. Notice that in this case the polar 
surface is not necessarily minimal.

\section{The results for surfaces}

In this section, we state our results on the associated family to an elliptic 
surface with circular ellipses of curvature.  We assert the existence of 
the associated family and discuss when the family is trivial. Then, we state
a general result that shows that the families associated to two consecutive 
circular ellipses coincide.

\begin{theorem}\po\label{main}
Let $g\colon L^2\to\Q_{c}^N$, $N\geq 6$,  be a simply connected 
substantial elliptic surface with $\mathcal{E}_\ell^g$ a circle for some 
$1\leq\ell\leq\tau^o-1$. Then there exists an associated one-parameter
family of elliptic surfaces with respect to the same almost complex structure
$$
G_\ell=\{g_\theta\colon L^2\to\Q_{c}^N:\theta\in\Sp^1=[0,\pi)\}
$$ 
with $\mathcal{E}_\ell^{g_\theta}$ a circle and such that
for each $\theta\in\Sp^1$ there exists a vector bundle isometry 
$\phi_\theta\colon N_gL\to N_{g_\theta}L$  that preserves the fundamental 
forms $\alpha^k_{g_\theta}=\phi_\theta\,\alpha^k_g$ for $2\leq k\leq\ell+1$ 
as well as the normal curvature tensor.  Moreover, any pair of elements 
in $G_\ell$ are non-congruent unless the family is trivial.
\end{theorem}

By the associate family $G_\ell$ being \emph{trivial} we mean that it only 
contains one element, that is, any $g_\theta$ is congruent to $g$ in the ambient 
space. The precise conditions for the family to be trivial are given by the following 
result. Notice that for $\ell=0$ the above associated family corresponds to the 
standard associated family to a minimal surface.

\begin{theorem}\po\label{trivial}
Let  $g\colon L^2\to\Q_{c}^N$, $N\geq 6$,  be a simply connected substantial 
elliptic surface such that $\mathcal{E}_\ell^g$ is a circle for some
$0\leq\ell\leq\tau^o-1$. If a pair of surfaces $g_{\theta},g_{\tilde{\theta}}\in G_\ell$ 
for $\theta\neq\tilde{\theta}$ are congruent, then $N$ is even and the  
$\mathcal{E}_s^g$'s are  circles for  $\ell\leq s\leq\tau$. 
Conversely, if $N$ be even and  the $\mathcal{E}_s^g$'s are 
circles for $\ell\leq s\leq\tau$, then the associated family $G_\ell$ 
is trivial.
\end{theorem}

In odd codimension the one-parameter associated family is never trivial.  
For minimal surfaces in a special case this was already observed in \cite{Jo} and  \cite{Vla}.
\medskip

The following is a consequence of the above result and basic  properties of polar 
surfaces of elliptic surfaces.

\begin{corollary}\po\label{cor1} Let $g\colon L^2\to\Q_{c}^N$, $c=0,1$, be a 
simply connected substantial elliptic surface where $N$ is even and $\mathcal{E}_\ell^g$ 
is a circle
for some $0\leq\ell\leq\tau-1$. Then, the associated family  $G_\ell$ is trivial 
if and only if $N$ is even and $g$ is  (locally) a polar surface to an $m$-isotropic surface
for $m=\tau-\ell$.
\end{corollary}

Finally, we prove the following result which shows that if two ellipses of 
consecutive order are circles then the associated families coincide.

\begin{theorem}\po\label{relation}
Let $\colon L^2\to\Q_{c}^N$ be a simply connected elliptic surface such that
the ellipses $\mathcal{E}_\ell^g,  \mathcal{E}_{\ell+r}^g$  are circles for some  
$0\leq \ell <\ell+r$. Then, the following facts are equivalent: 
\begin{itemize}
\item[(i)] $G_\ell \cap G_{\ell+r} \neq\{ g \}$.
\item[(ii)] The ellipses $\mathcal{E}_j^g$ are circles for $\ell\leq j\leq\ell+r$.  
\item[(iii)] $G_\ell = G_{\ell+r}$.
\end{itemize}
\end{theorem}

For the case of minimal surfaces we thus have the following.

\begin{corollary}\po\label{minimal}
Let  $g\colon L^2\to\Q_{c}^N$, $N\geq 6$, be a simply connected substantial minimal
surface with $\mathcal{E}_\ell^g$ a circle for some $1\leq\ell\leq\tau^o-1$. 
Then  $G_\ell=G_0$ if and only if $g$ is $\ell$-isotropic.
\end{corollary}

\subsection{The compatibility equations}

A key  ingredient in the proofs are the basic equations from the classical 
Burstin-Mayer-Allendoerfer theory discussed in Vol.\ IV of \cite{Sp}.  
They naturally extend the situation for curves under similar regularity 
conditions. The main result is that for a regular submanifold of a space form 
the tensors determined by the Frenet  equations  are a complete 
set of invariants.
\vspace{1,5ex} 

The Frenet equations for a regular \ii $f\colon M^n\to\Q_c^N$ are given by
$$
\nab_X\xi=-A^s_\xi X+D^s_X\xi+\mathsf{S}^s_X\xi\;\;
\mbox{if}\;\;\xi\in N^f_s\;\mbox{and}\;X\in T_xM,\;\;s\geq 1,
$$
in term of the linear maps
$$
A^s\colon TM\times N^f_s\to N^f_{s-1}\;\;\;\;\mbox{defined by}\;\;\;\; 
A^s_\xi X=-\pi_{s-1}(\nab_X\xi),
$$
$$
\hspace*{-4ex} D^s\colon TM\times N^f_s\to N^f_s
\;\;\;\;\;\;\;\mbox{defined by}\;\;\;\; 
D^s_X\xi=\pi_s(\nap_X\xi),
$$
$$
\hspace*{-2ex}\mathsf{S}^s\colon TM\times N^f_s\to N^f_{s+1}\;\;\;\;\;
\mbox{defined by}\;\;\;\;\mathsf{S}^s_X\xi=\pi_{s+1}(\nap_X\xi), 
$$
where $\nab$ is the connection in the induced bundle 
$ f^*(T\Q_c^N)=N^f_0\oplus N_fM$   and $\pi_{0}$ is the projection onto 
$N^f_0=f_*(TM)$. Notice that $A^1_{\xi}$ is the standard Weingarten operator
and that $D^s$ is a connection in $N^f_s$ compatible with the metric. 
An important fact is that  the tensors $A^s$ and $\mathsf{S}^s$ are 
completely determined by the higher fundamental forms since 
$$
\mathsf{S}^s_X(\alpha^{s+1}_f(X_1,\ldots,X_{s+1}))=\alpha^{s+2}_f(X,X_1,
\ldots,X_{s+1})
$$
and
\be\label{second}
\<A^s_\xi X,\eta\>=\<\xi,\mathsf{S}^{s-1}_X\eta\>\;\;
\mbox{for}\;\;\xi\in N^f_s\;\;\mbox{and}\;\;\eta\in N^f_{s-1}.
\ee

We  briefly summarize the basic results of the theory:
Let $f,\tilde{f}\colon M^n\to\Q_c^N$ be two regular isometric immersions. 
If there are vector bundle isometries $\phi_k\colon N_k^f\to N_k^{\tilde{f}}$ 
for all $k\geq 1$, which preserve the fundamental forms $\a^{k+1}$ and the induced 
normal connections $D^k$, then there is an isometry $\tau$ of $\Q_c^N$ such that 
$\tilde{f}=\tau\circ f$ and $\phi_k=\tau_*|_{N_k^f}$. Moreover, there is a set of 
equations given below, namely, the Generalized Gauss and Codazzi equations, 
that relate the higher fundamental forms and the induced connections.  
It turns our that the set of connections $D^k$ in $N^f_k$ is the unique set 
for which the higher order fundamental forms satisfy the Codazzi equations. 
Furthermore, the Generalized Gauss and Codazzi equations are the integrability conditions 
that assure the existence of an \ii provided all data involved has been provided.
\vspace{1,5ex}

\noindent The \underline{Generalized Gauss equation}.
\be\label{gengauss}
A^{s+1}_{\mathsf{S}^s_Y\xi}X-A^{s+1}_{\mathsf{S}^s_X\xi}Y
=D^s_XD^s_Y\xi-D^s_YD^s_X\xi-\mathsf{S}^{s-1}_XA^s_\xi Y 
+\mathsf{S}^{s-1}_YA^s_\xi X-D^s_{[X,Y]}\xi
\ee
for all $X,Y\in TM$ and $\xi\in N^f_s$. 
\medskip

\noindent
The \underline{Generalized Codazzi equation}. 
\be\label{gencodazzi}
D^{s+1}_X(\mathsf{S}^s_Y\xi)-D^{s+1}_Y(\mathsf{S}^s_X\xi)
+\mathsf{S}^s_XD^s_Y\xi-\mathsf{S}^s_YD^s_X \xi-\mathsf{S}^s_{[X,Y]}\xi=0
\ee
for all $X,Y\in TM$ and $\xi\in N^f_s$. \vspace{1,5ex}

\noindent Using  (\ref{second}) we have that (\ref{gencodazzi}) has the equivalent form 
\be\label{gencodazzi2}
D^s_XA^{s+1}_\xi Y-D^s_YA^{s+1}_\xi X 
+ A^{s+1}_{D^{s+1}_Y\xi}X - A^{s+1}_{D^{s+1}_X\xi}Y-A^{s+1}_{\xi}[X,Y]=0
\ee
for all $X,Y\in TM$ and $\xi\in N^f_{s+1}$.\vspace{1,5ex}

We conclude with some useful symmetric equations.

\begin{proposition}\po It holds that
\be\label{sym}
\mathsf{S}_Y^{s+1}\mathsf{S}_X^s\xi=\mathsf{S}_X^{s+1}\mathsf{S}_Y^s\xi\;\;
\mbox{or, equivalently, that}\;\;A_{A_\xi^sX}^{s-1}Y=A_{A_\xi^sY}^{s-1}X
\ee
for any $\xi\in N^f_s$ and $X,Y\in TM$.
\end{proposition}

\proof To prove the first equation take $\xi=\alpha^{s+1}_f(X_1,\ldots,X_{s+1})$
and use the symmetry of the higher fundamental forms. For the proof of the equivalent 
second equation take $\xi\in N^f_s$, $\eta\in N^f_{s-2}$ and use (\ref{second}) 
twice to obtain
$$
\<A^{s-1}_{A^s_\xi X}Y-A^{s-1}_{A^s_\xi Y}X,\eta\>
=\<\xi,\mathsf{S}^{s-1}_X\mathsf{S}^{s-2}_Y\eta
-\mathsf{S}^{s-1}_Y\mathsf{S}^{s-2}_X\eta\>=0,
$$
and this concludes the proof.\qed

\subsection{The proofs}

For a substantial elliptic surface $g\colon L^2\to\Q_{c}^N$ with a circular 
ellipse of curvature in a space form we first define a one-parameter 
family of compatible connections. Hereafter, we assume that $\mathcal{E}_\ell^g$
is a circle for given $0\leq\ell\leq\tau^o-1$, that is, the almost 
complex structure $J_\ell$ ($J_0=J$) is a vector bundle isometry. 
Notice that $J_\ell$ is  parallel with respect to the induced 
connection on $N^g_\ell$ by dimension reasons. Thus, for any 
$\varphi\in\Sp^1=[0,\pi)$ the map 
$R^\ell_\varphi\colon N^g_\ell\to N^g_\ell$ 
defined by (\ref{rotation}) is also a parallel isometry, i.e., 
\be\label{D0}
\pi_\ell(\nap_XR^\ell_\varphi\xi)=R^\ell_\varphi\pi_\ell(\nap_X\xi).
\ee

Let $\nab^\theta\colon TL\times g^*T\Q_c^N\to g^*T\Q_c^N$ for each $\theta\in \Sp^1$
be the map defined by modifying the induced connection $\nab$ of $g^*T\Q_c^N$ as follows:  
$$
\left\{\begin{array}{l}
\pi_{\ell+1}(\nab^\theta_X\xi)
=\pi_{\ell+1}(\nab_XR^\ell_{\theta}\xi)
\;\;\mbox{if}\;\;\;\xi\in N^g_\ell
\vspace*{1.5ex}\\
\pi_\ell(\nab^\theta_X\eta)= R^\ell_{-\theta}\pi_\ell(\nab_X\eta)
\;\;\;\;\:\;\mbox{if}\;\;\;\eta\in N^g_{\ell+1},
\end{array} \right.\nonumber
$$
and $\nab^\theta=\nab$ in all other cases. 
We also define the map $\nat\colon TL\times N_gL\to N_gL$ by 
$$
\nat_X \xi= \nab^\theta_X \xi-\pi_0 (\nab^\theta_X \xi).
$$

For $\ell=0$, we have the map $\a_\theta\colon TL\times TL\to N_gL$ given by
$$
a_\theta (X,Y)=\pi_1(\nab^\theta_X g_*Y).
$$
Thus $a_\theta (X,Y)=\a_g(J_\theta X, Y)$ where $J_\theta=\cos\theta I+\sin \theta J$.
Also $\nat=\nap$. Then,  the triple 
$(\alpha_\theta, \<\,,\,\>,\nap)$ satisfies the Gauss, 
Codazzi and Ricci equations. Therefore, if  $L$ is simply connected if follows from 
the Fundamental theorem of submanifolds that there exists an isometric minimal surface 
$g_\theta\colon L\to\Q_c^N$ and a parallel vector bundle isometry
$$
\phi_\theta\colon (N_gL,\nat)\to (N{g_\theta}L,\nabla^\perp(g_\theta))
$$
such that $\alpha_{g_\theta}=\phi_\theta\,\alpha_\theta$.
Of course, the family $g_\theta$ with $\theta\in\Sf^1$ is just the standard 
associated family of the minimal surface $g$.

For $\ell\geq 1$, the map $\nat$ is 
obtained modifying the normal connection of $g$ as follows:  
\be\label{cone}
\left\{\begin{array}{l}
\pi_{\ell+1}(\nat_X\xi)
=\pi_{\ell+1}(\nap_XR^\ell_{\theta}\xi)
\;\;\mbox{if}\;\;\;\xi\in N^g_\ell
\vspace*{1.5ex}\\
\pi_\ell(\nat_X\eta)= R^\ell_{-\theta}\pi_\ell(\nap_X\eta)
\;\;\;\;\:\;\mbox{if}\;\;\;\eta\in N^g_{\ell+1},
\end{array} \right.
\ee
and $\nat=\nap$ in all other cases. 

\begin{lemma}\po\label{ral} For $\ell \geq 1$ the map $\nat$ is 
a Riemannian connection whose curvature tensor satisfies \mbox{$R^\theta=R^\perp$.}
\end{lemma}

\proof  Take  $\xi\in N^g_\ell$ and $\eta\in N^g_{\ell+1}$. Then,
\be\label{x}
\nat_Xf\xi=
(\pi_{\ell-1}+\pi_\ell)(\nap_Xf\xi)+\pi_{\ell+1}(\nap_XfR^\ell_\theta\xi)
=X(f)\xi+f\nat_X \xi
\ee
and
\be\label{xx}
\nat_Xf\eta=(\pi_{\ell+1} +\pi_{\ell+2})(\nap_Xf\eta) 
+R^\ell_{-\theta}\pi_\ell(\nap_Xf\eta)=X(f)\eta+f\nat_X\eta.
\ee
Moreover, we obtain using (\ref{D0}) that
\be\label{xxx}
\<\nat_X\xi,\eta\>+\<\xi,\nat_X\eta\>\!=\!\<\nap_XR^\ell_\theta\xi,\eta\>
+\<\xi,R^\ell_{-\theta}\nap_X\eta\>\!=\!\<R^\ell_\theta\nap_X\xi,\eta\>
+\<\xi,R^\ell_{-\theta}\nap_X\eta\>=0,
\ee
and that the connection is Riemannian follows easily 
from (\ref{x}), (\ref{xx}) and (\ref{xxx}).

The second claim amounts to show  that the tensor defined by
$$
\Ral(X,Y)\xi=R^\theta(X,Y)\xi-R^\perp (X,Y)\xi
$$
vanishes. In the sequel, some of the arguments will just be sketched
to avoid writing rather long but straightforward computations.
\medskip 

\noindent We divide the proof in several cases:
\medskip 

\noindent\emph{Case 1.} The case $\xi\in N^g_1
\oplus\cdots\oplus N^g_{\ell-2}\oplus N^g_{\ell+3}\oplus\cdots\oplus N^g_{\tau}$
is trivial.
\medskip 

\noindent\emph{Case 2.} Take $\xi\in N^g_{\ell+2}$. Then,
$$
\Ral(X,Y)\xi=
(R^{\ell}_{-\theta}-I)\Big(A_{A_\xi^{\ell+2}Y}^{\ell+1}X
-A_{A_\xi^{\ell+2}X}^{\ell+1}Y\Big),
$$
and the claim follows from (\ref{sym}).
\medskip

\noindent\emph{Case 3.} Take $\xi\in N^g_{\ell+1}$. Then,
$$
\Ral(X,Y)\xi=B_\theta(X,Y)-B_\theta(Y,X)-B_0(X,Y)+B_0(Y,X)
+(I-R^\ell_{-\theta})\pi_\ell(\nabla_{[X,Y]}^\perp \xi)
$$
where 
$$
B_\theta(X,Y)=
R^\ell_{-\theta}\pi_\ell(\nap_X \pi_{\ell+1}(\nap_Y\xi))
+\pi_\ell(\nap_XR^\ell_{-\theta}\pi_\ell(\nap_Y\xi))+
\pi_{\ell-1}(\nap_XR^\ell_{-\theta}\pi_\ell(\nap_Y\xi)).
$$
Observe that  (\ref{D0}) can be written as
\be\label{D}
D^\ell_XR^\ell_\varphi\xi=R^\ell_\varphi D^\ell_X\xi.
\ee
We obtain using (\ref{sym}) and (\ref{D})  that 
\bea
\Ral(X,Y)\xi\!\!\!&=&\!\!\!(I-R^\ell_{-\theta})\Big(D^{\ell}_XA^{\ell+1}_\xi Y
-D^{\ell}_YA^{\ell+1}_\xi X
+A^{\ell+1}_{D^{\ell+1}_Y \xi}X
-A^{\ell+1}_{D^{\ell+1}_X\xi}Y
-A^{\ell+1}_\xi {[X,Y]}\Big)\\
\!\!\!&&\!\!\!
+A^\ell_{R^\ell_{-\theta}A^{\ell+1}_{\xi}Y} X
-A^\ell_{R^\ell_{-\theta}A^{\ell+1}_{\xi}X} Y.
\eea
For $\eta\in N^g_{\ell-1}$, we have using (\ref{second}) that 
\be\label{as}
\<A^\ell_{R^\ell_{-\theta}A^{\ell+1}_{\xi}Y}X,\eta\> 
=\<A^{\ell+1}_\xi Y,{R^\ell_{\theta}\mathsf{S}^{\ell-1}_X\eta\>
=\<\xi,}\mathsf{S}^\ell_YR^\ell_{\theta}\mathsf{S}^{\ell-1}_X\eta\>.
\ee
Since 
$R^\ell_\theta\mathsf{S}^{\ell-1}_X\xi
=\mathsf{S}^{\ell-1}_XR_\theta^{\ell-1}\xi$ from (\ref{one0}),
we obtain from (\ref{sym}) that
\be\label{srs}
\mathsf{S}^\ell_Y R^\ell_{\theta}\mathsf{S}^{\ell-1}_X\xi
=\mathsf{S}^\ell_X R^\ell_{\theta}\mathsf{S}^{\ell-1}_Y\xi.
\ee 
Now the claim follows from (\ref{gencodazzi2}), (\ref{as}) and (\ref{srs}).
\medskip

\noindent\emph{Case 4.} Take $\xi\in N^g_{\ell}$. 
First assume $\ell\geq 2$. Using (\ref{gengauss}), 
(\ref{gencodazzi}),  (\ref{sym}) and (\ref{D}) we obtain
\be\label{z}
R^{\theta}(X,Y)\xi 
=\mathsf{S}^{\ell-1}_YA^\ell_\xi X
-\mathsf{S}^{\ell-1}_XA^\ell_\xi Y
- R^\ell_{-\theta}(\mathsf{S}^{\ell-1}_YA^\ell_{R^\ell_\theta\xi}X
-\mathsf{S}^{\ell-1}_XA^\ell_{R^\ell_\theta\xi}Y).
\ee
On the other hand, it holds that
\be\label{rsa}
R^\ell_\varphi(\mathsf{S}^{\ell-1}_YA^\ell_{\xi}X-\mathsf{S}^{\ell-1}_X
A^\ell_{\xi}Y)=\mathsf{S}^{\ell-1}_YA^\ell_{R^\ell_\varphi\xi}X
-\mathsf{S}^{\ell-1}_XA^\ell_{R^\ell_\varphi\xi}Y.  
\ee
In fact, it follows using (\ref{second}) that
$$
\<R^\ell_\varphi(\mathsf{S}^{\ell-1}_YA^\ell_{\xi}X-\mathsf{S}^{\ell-1}_X
A^\ell_{\xi}Y), R^\ell_\varphi\delta\>
=\<A^\ell_{\delta}Y,A^\ell_{\xi}X\>-\<A^\ell_{\delta}X,A^\ell_{\xi}Y\>
$$
and
$$
\<\mathsf{S}^{\ell-1}_YA^\ell_{R^\ell_\varphi\xi}X
-\mathsf{S}^{\ell-1}_XA^\ell_{R^\ell_\varphi\xi}Y,R^\ell_\varphi\delta\>
 =\< A^\ell_{R^\ell_\varphi\delta}Y,A^\ell_{R^\ell_\varphi\xi}X\>
-\<A^\ell_{R^\ell_\varphi\delta}X,A^\ell_{R^\ell_\varphi\xi}Y\>.
$$
Since $N^g_\ell=\spa\{\xi,J_\ell\xi\}$, to obtain
(\ref{rsa}) is suffices to compute  the right hand side of 
both equations for $\delta=J_\ell\xi$ and observe that they coincide.

We now have from (\ref{z}) and (\ref{rsa}) that $R^{\theta}(X,Y)\xi=0$, 
and this proves the claim since also $R^{\perp}(X,Y)\xi=0$ from the 
Ricci equation.
\medskip

\noindent For  $\ell=1$, we have that
$$
R^{\theta}(X,Y)\xi =   
R^\ell_{-\theta}\left(\alpha_g(X,A_{R^\ell_\theta\xi}Y)
-\alpha_g(Y,A_{R^\ell_\theta\xi}X)\right).
$$
Since $A_{R^\ell_\theta\xi}= R_\theta A_\xi$, we obtain 
from the Ricci equation that
$R^{\theta}(X,Y)\xi = R^{\perp}(X,Y)\xi$ 
and the claim also follows in this case. 
\medskip

\noindent\emph{Case 5.} Take $\xi\in N^g_{\ell-1}$.  Then,  
$$
\Ral(X,Y)\xi=B_\theta(X,Y)-B_\theta(Y,X)-B_0(X,Y)+B_0(Y,X)
$$
where
$$
B_\theta(X,Y)=\pi_\ell(\nap_X \pi_{\ell-1}(\nap_Y\xi))
+\pi_\ell(\nap_X\pi_\ell(\nap_Y\xi))+
\pi_{\ell+1}(\nap_XR^\ell_\theta\pi_\ell(\nap_Y\xi)).
$$
It follows that 
\bea
\Ral(X,Y)\xi\!\!\!&=&\!\!\!(R^\ell_{-\theta}-I)
\left(D^{\ell}_X(\mathsf{S}^{\ell-1}_Y\xi)
-D^{\ell}_Y(\mathsf{S}^{\ell-1}_X\xi)
+\mathsf{S}^{\ell-1}_XD^{\ell-1}_Y\xi-\mathsf{S}^{\ell-1}_YD^{\ell-1}_X \xi
-\mathsf{S}^{\ell-1}_{[X,Y]}\xi\right)\\
\!\!\!&&\!\!\!+\mathsf{S}^{\ell}_X R^\ell_\theta\mathsf{S}^{\ell-1}_Y\xi
-\mathsf{S}^{\ell}_Y R^\ell_\theta\mathsf{S}^{\ell-1}_X\xi,
\eea 
and the claim follows from (\ref{gencodazzi}) and (\ref{srs}).
\medskip

To conclude the proof, we observe that the case $\tau^o=\tau-1$
and $\xi\in N^g_{\tau}$ is included in the above cases.\vspace{1,5ex}\qed

\noindent{\it Proof of Theorem \ref{main}.} 
We argue that the triple $(\alpha_g, \<\,,\,\>,\nat)$ satisfies the Gauss, 
Codazzi and Ricci equations.  Then, according to the Fundamental theorem of 
submanifolds there exists an isometric immersion 
$g_\theta\colon L^2\to\Q_{c}^N$ and a parallel vector bundle isometry
\be\label{fts}
\phi_\theta\colon (N_gL,\nat)\to (N{g_\theta}L,\nabla^\perp(g_\theta))
\ee
such that $\alpha_{g_\theta}=\phi_\theta\,\alpha_g$.

The Gauss equation holds since the second fundamental
form remains the same.  The Codazzi equation
$$
(\nabla_X^{\theta}\alpha_g)(Y,Z)=(\nabla_Y^{\theta}\alpha_g)(X,Z)
$$
is trivially satisfied for $\ell\geq 3$ since $\nabla^\theta\alpha_g
=\nabla^\perp\alpha_g$. For $\ell=1$, we obtain
$$ 
(\nabla_X^{\theta}\alpha_g)(Y,Z)=\pi_1((\nabla_X^\perp\alpha_g)(Y,Z))
+\pi_2((\nabla_X^\perp\alpha_g)(Y,R_\theta Z))
$$
while for $\ell=2$, we have 
$$ 
(\nabla_X^{\theta}\alpha_g)(Y,Z)=(\pi_1+\pi_2)((\nabla_X^\perp\alpha_g)(Y,Z)),
$$
and again the Codazzi equation follows.

The proof that the curvature tensor $R^{\theta}$ of $\nat$ satisfies the 
Ricci equation
$$
R^{\theta}(X,Y)\xi=\alpha_g(X,A_{\xi}Y)-\alpha_g(Y,A_{\xi}X)
$$
is  a consequence of Lemma \ref{ral} above.
Finally, the statements on the fundamental forms is part of Proposition \ref{HFF}
given next.\vspace{1,5ex}\qed

The following result provides the expressions for the higher fundamental 
forms  of the associated family $g_\theta$ in terms of the ones 
corresponding   to $g$. We observe that this can be used to give an alternative 
(but more complicate) definition for the associated family by means of the 
version in \cite{Sp} of the Fundamental theorem of submanifolds as part of 
the Burstin-Mayer-Allendoerfer theory.

\begin{proposition}\po\label{HFF} Let  $g\colon
L^2\to\Q_{c}^N$
be a simply connected elliptic surface with $\mathcal{E}_\ell^g$
a circle for some $1\leq\ell\leq\tau^o-1$.
Then, up to identification, the higher fundamental forms of 
$g_\theta$ are given by
$$
\alpha_{g_\theta}^s(X_1,\ldots,X_s)\!=\!\left\{\begin{array}{l}
\alpha_g^s(X_1,\ldots,X_s)\;\;\mbox{if}\;\;
2\leq s\leq\ell+1,
\vspace*{1.5ex}\\
R^{s-1}_\theta\alpha_g^s(X_1,\ldots,X_s)
=\alpha_g^s(R_\theta X_1,\ldots,X_s)
\;\;\mbox{if}\;\;\ell+2\leq s\leq\tau^o+1.
\end{array} \right.
$$
Moreover, if $N$ is odd, then 
$$
\alpha_{g_\theta}^{\tau+1}(X_1,\ldots,X_{\tau+1})=\alpha_g^{\tau+1}(R_\theta X_1,\ldots,X_{\tau+1}).
$$
\end{proposition}

\proof  Since $\alpha^2_{g_\theta}=\alpha_g$, the case $2\leq s\leq\ell$ 
follows easily from the definitions. For the other cases, we have 
\bea
\alpha^{\ell+1}_{g_\theta}(X_1,\ldots,X_{\ell+1})
\!\!\!&=&\!\!\!\pi_\ell (\nat_{X_{\ell+1}} \alpha^\ell
_{g_\theta}(X_1,\ldots, X_\ell))=\pi_{\ell}(\nat_{X_{\ell+1}}
\alpha^{\ell}_{g}(X_1,\ldots,X_\ell))\\
\!\!\!&=&\!\!\!\pi_\ell (\nap_{X_{\ell+1}}
\alpha^{\ell}_{g}(X_1,\ldots,X_\ell))
=\alpha^{\ell+1}_g(X_1,\ldots,X_{\ell+1})
\eea
and 
\bea
\alpha^{\ell+2}_{g_\theta}(X_1,\ldots,X_{\ell+2})
\!\!\!&=&\!\!\! \pi_{\ell+1} (\nat_{X_{\ell+2}} \alpha^{\ell+1}
_{g_\theta}(X_1,\ldots, X_{\ell+1}))\\
\!\!\!&=&\!\!\!\pi_{\ell+1}(\nat_{X_{\ell+2}}
\alpha^{\ell+1}_{g}(X_1,\ldots,X_{\ell+1}))\\
\!\!\!&=&\!\!\!\pi_{\ell+1} (\nap_{X_{\ell+2}}R^\ell_\theta
\alpha^{\ell+1}_{g}(X_1,\ldots,X_{\ell+1}))\\
\!\!\!&=&\!\!\!\pi_{\ell+1}(\nap_{X_{\ell+2}}\alpha^{\ell+1}_g
(R_\theta X_1,\ldots,X_{\ell+1}))\\
\!\!\!&=&\!\!\!\alpha^{\ell+2}_{g}(R_\theta X_1,\ldots,X_{\ell+2})\\
\!\!\!&=&\!\!\!R^{\ell+1}_\theta\alpha^{\ell+2}_{g}(X_1,\ldots,X_{\ell+2}),
\eea
and the remaining of the proof is immediate.\vspace{1,5ex}\qed

\noindent{\it Proof of Theorem \ref{trivial}.}  Since the case of $\ell=0$ is 
well-known we argue for $\ell\geq 1$.
Suppose first that the surfaces $g_{\theta}$ and $g_{\tilde{\theta}}$ in $G_\ell$ 
are congruent.  Without loss of generality, we may assume that $\tilde{\theta}=0$.
Then, there exists a parallel vector bundle isometry
$\psi\colon N_gL\to N_{g_\theta}L$ such that $\psi(N^g_s)=N^{g_\theta}_s$ and 
$\alpha^s_{g_\theta}=\psi\,\alpha^s_{g}$ for any $2\leq s\leq\tau^o$. 
Proposition \ref{HFF} yields 
$$
\psi\,\alpha_{g}^s=\left\{\begin{array}{l}
\alpha_g^s\;\;\;\;\;\;\;\;\;\;\mbox{if}\;\;2\leq s\leq\ell+1,
\vspace*{1.5ex}\\
R^{s-1}_\theta\,\alpha_g^s
\;\;\mbox{if}\;\;\ell+2\leq s\leq\tau^o+1,
\end{array} \right.
$$
and if $N$ is odd then
$$
\psi\alpha_{g}^{\tau+1}(X_1,\ldots,X_{\tau+1})
=\alpha_g^{\tau+1}(R_\theta X_1,\ldots,X_{\tau+1}).
$$
Therefore, we have that $N$ is even and 
\be\label{def}
\psi= I\;\mbox{on}\; N^g_1\oplus\cdots\oplus N^g_\ell\;\;\mbox{and}\;\; 
\psi=R^s_\theta\;\mbox{on}\; N^g_s\; \;\mbox{if}\;s\geq\ell+1.
\ee
We conclude that $R^s_\theta$ is an isometry for
$s\geq\ell+1$ and hence $J_s$ is an isometry for $s\geq\ell$.
\vspace{1ex}

Conversely, suppose that $N$ is even and any $\mathcal{E}_s(g)$ is a circle for all
$s\geq\ell$, or equivalently, that $R^s_\theta$ is an isometry for
$s\geq\ell$. Then $\psi\colon N_gL\to N_{g_\theta}L$ given by (\ref{def})
is a vector bundle isometry that preserves the second fundamental form. 
To conclude the proof it remains to show that $\psi$ is parallel, i.e., 
$\psi \nap_X\xi=\nat_X\psi\xi$.  For that we distinguish several cases:
\medskip

\noindent\emph{Case 1.} The case  $\xi\in N^g_1\oplus\ldots\oplus N^g_{\ell-1}$
is trivial.
\medskip

\noindent\emph{Case 2.}  Assume that $\xi\in N^g_{\ell}$. 
We have using (\ref{one0}) that
\bea
\nat_X\psi\xi \!\!\!&=&\!\!\! \nat_X\xi
=(\pi_{\ell-1}+ \pi_\ell)(\nap_X\xi)
+\pi_{\ell+1}(\nap_X R^{\ell}_\theta\xi)\\
\!\!\!&=&\!\!\!(\pi_{\ell-1}+ \pi_\ell)(\nap_X\xi)
+ R^{\ell+1}_\theta\pi_{\ell+1}(\nap_X\xi)=\psi \nap_X \xi.
\eea

\noindent\emph{Case 3.} Assume that $\xi\in N^g_{\ell+1}$. We have, 
$$
\nat_X\psi\xi=\nat_XR^{\ell+1}_\theta\xi=
R^\ell_{-\theta}\pi_\ell(\nap_XR^{\ell+1}_\theta\xi)
+(\pi_{\ell+1}+\pi_{\ell+2})(\nap_XR^{\ell+1}_\theta\xi)
$$
and
$$
\psi\nap_X\xi=\pi_\ell(\nap_X\xi)
+R^{\ell+1}_\theta\pi_{\ell+1}(\nap_X\xi)
+ R^{\ell+2}_\theta\pi_{\ell+2}(\nap_X\xi).
$$
To obtain equality we observe that (\ref{two0}) yields
$$
R^\ell_{-\theta}\pi_\ell(\nap_XR^{\ell+1}_\theta\xi)=\pi_\ell(\nap_X\xi)
$$
and that the $N^g_{\ell+1}$  and  $N^g_{\ell+2}$ components are equal 
due to (\ref{D0}) and (\ref{one0}), respectively.\vspace{1ex}

\noindent\emph{Case 4.} Assume that $\xi\in N^g_s$ for $s\geq\ell+2$.
We have,
$$
\nat_X\psi\xi=\nap_XR^s_\theta\xi
=(\pi_{s-1}+\pi_s+\pi_{s+1})(\nap_XR^s_\theta\xi),
$$
$$
\psi\nap_X\xi=R^{s-1}_\theta\pi_{s-1}(\nap_X\xi)
+R^s_\theta\pi_s(\nap_X\xi)
+R^{s+1}_\theta\pi_{s+1}(\nap_X\xi),
$$
and equality follows from (\ref{one0}), (\ref{two0}) and (\ref{D0}).
\vspace{1,5ex}\qed

\noindent{\it Proof of Corollary \ref{cor1}.} The proof follows from 
Theorem \ref{trivial}  and the fact proved in \cite{DF2} that
an elliptic surface has circular curvature ellipses from some order on if
and only if any polar surface has circular curvature ellipses up to that order.
\qed\vspace{1,5ex}

\noindent{\it Proof of Theorem \ref{relation}.}  
To see that  $(i)\Rightarrow (ii)$ suppose that 
a surface $g_\theta$ in $G_\ell$ is congruent to $\tilde{g}_\omega$ 
in $G_{\ell+r}$.  First assume that $\ell\geq 1$, and let 
$$
\phi_\theta\colon (N_gL,\nabla^{\theta,\ell}\to (N_{g_\theta}L,\nabla^\perp(g_\theta)),\;\;\;
\sigma_\omega\colon (N_gL,\nabla^{\omega,\ell+r})
\to (N_{\tilde{g}_\omega}L,\nabla^\perp(\tilde{g}_\omega))
$$
be the parallel vector bundle  isometries given by (\ref{fts}). 
By assumption, there exists a  parallel bundle  isometry 
$$
\psi\colon (N_{g_\theta}L,\nabla^\perp(g_\theta))
\to (N_{\tilde{g}_\omega}L,\nabla^\perp(\tilde{g}_\omega)),
$$
such that
$\a^{s+1}_{\tilde{g}_\omega}=\psi\, \a^{s+1}_{g_\theta}$
for any $1\leq s\leq \tau-1$. Proposition \ref{HFF} yields
$$
\a^{s+1}_{\tilde{g}_\omega}=\sigma_\omega\,\a^{s+1}_g
$$
for any $1\leq s\leq \ell+r$ and 
$$
\a^{s+1}_{g_\theta}=\phi_\theta\, R^s_\theta\,\a^{s+1}_g
$$
for any $\ell+1\leq s\leq\ell+r$. 
We obtain that $\sigma_\omega\,\a^{s+1}_g
=\psi\circ\phi_\theta\circ R^s_\theta\,\a^{s+1}_g$, that is,
$$
R^s_\theta=(\psi\circ\phi_\theta)^{-1}\circ\sigma_\omega
$$ 
on $N^s_g$ for $\ell+1\leq s\leq\ell+r$.
Since $R^s_\theta$ is an isometry for any $\ell +1 \leq s\leq \ell+r$,
then all the $\mathcal{E}_j^g,\ell\leq j\leq\ell+r$,
are circles.

Now assume  $\ell=0$ and let  
$\phi_\theta\colon (N_gL,\nap)\to (N_{g_\theta}L,\nabla^\perp(g_\theta))$ 
be the parallel bundle isometry such that
$$
\a^{s+1}_{g_\theta}=\phi_\theta\circ R^s_{\theta}\,\a^{s+1}_g
$$
for any $1\leq s\leq \tau^o$. Proposition \ref{HFF} yields
$\a^{s+1}_{\tilde{g}_\omega}=\sigma_\omega\,\a^{s+1}_g$
for $1\leq s\leq \ell+r$. Then, 
$$
\a^{s+1}_{\tilde{g}_\omega}
=\sigma_\omega\circ (\phi_\theta\circ R^s_{\theta})^{-1} \a^{s+1}_{g_\theta}.
$$
Since
$\a^{s+1}_{\tilde{g}_\omega}=\psi\,\a^{s+1}_{g_\theta}$
we find that 
$R^s_\theta=(\psi\circ\phi_\theta)^{-1}\circ\sigma_\omega$ on $N^s_g$
for any $1\leq s\leq \ell+r$. Hence, all the  $\mathcal{E}_j^g, 0\leq j\leq\ell+r$,
are circles.
\bigskip

\noindent We show that $(ii) \Rightarrow (iii)$. 
At first, we assume $\ell=0$ and prove that $G_0=G_r$. According to Theorem \ref{main}, 
the second fundamental form of $g_\theta\in G_r$ is given by
$\alpha_{g_\theta}=\psi_\theta\alpha_g$.
The second fundamental form of $h_\theta \in G_0$ is given by 
$$
\a_{h_\theta}=\phi_\theta\circ R^1_\theta\,\a_g,
$$
where $\phi_\theta\colon (N_gL, \nap)\to (N_{h_\theta}L,\nap (h_\theta))$ is a parallel 
vector bundle isometry.  Hence $\a_{h_\theta} =T_{\theta}\, \a_{g_\theta}$
where $T_\theta\colon N_{g_\theta}L\to N_{h_\theta}L$ is the bundle isometry given by
$$
T_\theta=\phi_\theta\circ R^s_{\theta}\circ\psi_\theta ^{-1} \;
\mbox{on}\; N^{g_\theta}_s\;\;\mbox{if}\;1\leq s\leq r\;\;\mbox{and}\;\; 
T_\theta=\phi_\theta\circ\psi_\theta ^{-1}\;
\mbox{on}\; N^{g_\theta}_{r+1}\oplus\cdots\oplus N^{g_\theta}_{\tau}.
$$

We show next that $T_\theta$ is parallel, i.e.,
\be\label{p}
\nap_X T_\theta \xi_\theta=T_\theta(\nap_X \xi_\theta) 
\ee
where $\xi_\theta=\psi_\theta\,\xi$ and $\xi \in N^g_s, 1\leq s \leq \tau$. 
We need to distinguish several cases:
\medskip 

\noindent\emph{Case 1.}  Assume that $1\leq s\leq r-1.$ We have that
\bea
\nap_XT_\theta\xi_\theta\!\!\!&=&\!\!\!\nap_X\phi_\theta R^s_\theta\xi
=\phi_\theta(\nap_XR^s_\theta\xi).
\eea
Using (\ref{one0}), (\ref{two0}), (\ref{D0}) and (\ref{cone}) we obtain
\bea
T_\theta(\nap_X\xi_\theta)
\!\!\!&=&\!\!\!T_\theta\circ\psi_\theta(\nabla_X^{\theta,r}\xi)
=T_\theta\circ\psi_\theta(\pi_{s-1}(\nabla_X^{\theta,r}\xi)
+\pi_s(\nabla_X^{\theta,r}\xi)+\pi_{s+1}(\nabla_X^{\theta,r}\xi))\\
\!\!\!&=&\!\!\!T_\theta\circ\psi_\theta(\pi_{s-1}(\nap_X\xi)
+\pi_s(\nap_X\xi)+\pi_{s+1}(\nap_X\xi))\\
\!\!\!&=&\!\!\!\phi_\theta\circ R^{s-1}_\theta(\pi_{s-1}(\nap_X\xi))
+\phi_\theta \circ R^s_\theta(\pi_s(\nap_X\xi))+
\phi_\theta\circ R^{s+1}_\theta(\pi_{s+1}(\nap_X\xi))\\
\!\!\!&=&\!\!\!\phi_\theta((R^{s-1}_{-\theta})^t\pi_{s-1}(\nap_X\xi)
+\pi_s(\nap_X R^s_\theta\xi)
+\pi_{s+1}(\nap_XR^s_\theta\xi))\\
\!\!\!&=&\!\!\!\phi_\theta(\pi_{s-1}(\nap_XR^s_\theta\xi)
+\pi_s(\nap_XR^s_\theta\xi)
+\pi_{s+1}(\nap_XR^s_\theta \xi)),
\eea
and this proves (\ref{p}).

\medskip
\noindent\emph{Case 2.} Assume that $s=r$. As before, we have 
$$
\nap_X T_\theta\xi_\theta=\phi_\theta(\nap_X R^r_\theta\xi).
$$
Using  (\ref{one0}), (\ref{two0}), (\ref{D0}) and (\ref{cone}) 
it follows that
\bea
T_\theta(\nap_X\xi_\theta)\!\!\!&=&\!\!\! T_\theta\circ\psi_\theta(\nabla_X^{\theta, r}\xi)
=T_\theta\circ\psi_\theta(\pi_{r-1}(\nabla_X^{\theta,r}\xi)+\pi_r(\nabla_X^{\theta, r}\xi)
+\pi_{r+1}(\nabla_X^{\theta, r}\xi))\\
\!\!\!&=&\!\!\!T_\theta\circ\psi_\theta(\pi_{r-1}(\nap_X\xi)
+\pi_{r}(\nap_X\xi)+\pi_{r+1}(\nap_X R^r_\theta\xi))\\
\!\!\!&=&\!\!\!\phi_\theta\circ R^{r-1}_\theta (\pi_{r-1}(\nap_X\xi))
+\phi_\theta\circ R^r_\theta(\pi_r(\nap_X\xi))
+\phi_\theta(\pi_{r+1}(\nap_XR^r_\theta\xi))\\
\!\!\!&=&\!\!\!\phi_\theta((R^{r-1}_{-\theta})^t\pi_{r-1}(\nap_X\xi)
+\pi_r(\nap_XR^r_\theta\xi)+\pi_{r+1}(\nap_XR^r_\theta\xi))\\
\!\!\!&=&\!\!\!\phi_\theta(\pi_{r-1}(\nap_XR^r_\theta \xi)
+\pi_{r}(\nap_X R^{r}_\theta  \xi)+
\pi_{r+1}(\nap_XR^r_\theta \xi)),
\eea
and this  proves (\ref{p}).

\medskip
\noindent\emph{Case 3.} Assume that $s=r+1$. We have that 
$$
\nap_X T_\theta\xi_\theta=\phi_\theta(\nap_X\xi).
$$
On the other hand,
\bea
T_\theta(\nap_X\xi_\theta)
\!\!\!&=&\!\!\! T_\theta\circ\psi_\theta(\nabla_X^{\theta,r}\xi)
=T_\theta\circ\psi_\theta(\pi_r(\nabla_X^{\theta,r}\xi)
+\pi_{r+1}(\nabla_X^{\theta,r}\xi)
+\pi_{r+2}(\nabla_X^{\theta,r}\xi))\\
\!\!\!&=&\!\!\!T_\theta\circ\psi_\theta(R^r_{-\theta}\pi_{r}(\nap_X\xi)
+\pi_{r+1}(\nap_X\xi)+\pi_{r+2}(\nap_X\xi))\\
\!\!\!&=&\!\!\! \phi_\theta\circ R^r_\theta\circ R^r_{-\theta}(\pi_r(\nap_X\xi))
+\phi_\theta(\pi_{r+1}(\nap_X\xi))+
\phi_\theta(\pi_{r+2}(\nap_X\xi))\\
\!\!\!&=&\!\!\!\phi_\theta(\pi_r(\nap_X\xi)
+\pi_{r+1}(\nap_X\xi)+\pi_{r+2}(\nap_X\xi)),
\eea
and (\ref{p}) holds true.

\medskip
\noindent\emph{Case 4.} Assume that $s\geq r+2$. Then, we have
$$
\nap_XT_\theta\xi_\theta=\phi_\theta(\nap_X\xi)
$$
and 
\bea
T_\theta(\nap_X\xi_\theta)
\!\!\!&=&\!\!\!T_\theta\circ\psi_\theta(\nabla_X^{\theta,r}\xi)
=T_\theta\circ\psi_\theta(\pi_s(\nabla_X^{\theta,r}\xi)
+\pi_{s+1}(\nabla_X^{\theta,r}\xi)+\pi_{s+2}(\nabla_X^{\theta,r}\xi))\\
\!\!\!&=&\!\!\!T_\theta\circ\psi_\theta(\pi_{s}(\nap_X\xi)
+\pi_{s+1}(\nap_X\xi)+\pi_{s+2}(\nap_X\xi))\\
\!\!\!&=&\!\!\!\phi_\theta(\pi_s(\nap_X\xi)+\pi_{s+1}(\nap_X\xi)
+\pi_{s+2}(\nap_X\xi)).
\eea
Thus $T_\theta$ is parallel and $g_\theta\in G_0$. Hence  $G_r=G_0$. 
\medskip

Now assume  $\ell\geq 1$. 
Take $g_\theta\in G_\ell,h_\theta\in G_{\ell+r}$ and let 
$$
\phi_\theta\colon (N_gL,\nabla^{\theta, \ell})\to (N_{g_\theta},\nap (g_\theta)),\;\;\;
\psi_\theta\colon (N_gL,\nabla^{\theta, \ell+r})\to (N_{h_\theta},\nap (h_\theta))
$$ 
be the corresponding  parallel isometries given by  (\ref{fts}). 
Then, we have
$\a_{h_\theta}=S\,\a_{g_\theta}$ where $S\colon (N_{g_\theta}, 
\nap (g_\theta))\to(N_{h_\theta},\nap (h_\theta))$ is the bundle isometry 
given by
$$
S=\psi_\theta\circ\phi_\theta^{-1}\;\mbox{on}\; N^{g_\theta}_s\;\;
\mbox{if}\;1\leq s\leq\ell\;\;\text{or}\;\; s\geq\ell+r+1
$$
and
$$
S=\psi_\theta\circ R^s_{-\theta}\circ\phi_\theta ^{-1}\;\;
\mbox{on}\; N^{g_\theta}_s\;\;\mbox{if}\;\ell+1\leq s\leq\ell +r.
$$

We show next that $S$ is parallel, i.e., 
\be\label{pp}
\nap_XS\xi_\theta=S(\nap_X \xi_\theta)
\ee
where $\xi_\theta=\phi_\theta\,\xi$ and $\xi\in N^g_s, 1\leq s\leq\tau$. 
We need to distinguish several cases:
\medskip 

\noindent\emph{Case 1.}  Assume that $1\leq s \leq \ell-1.$ Then, we have
\be\label{s}
\nap_X S\xi_\theta=\nap_X\psi_\theta\xi
=\psi_\theta(\nabla_X ^{ \theta, \ell+r}\xi)
=\psi_\theta(\nap_X\xi).
\ee
We obtain
\bea
S(\nap_X \xi_\theta)
\!\!\!&=&\!\!\!S\circ\phi_\theta(\nabla_X^{\theta,\ell}\xi)
=S\circ\phi_\theta(\pi_{s-1}(\nabla_X^{\theta,\ell}\xi)
+\pi_s(\nabla_X^{\theta,\ell}\xi)+\pi_{s+1}(\nabla_X^{\theta,\ell}\xi))\\
\!\!\!&=&\!\!\!S\circ\phi_\theta(\pi_{s-1}(\nap_X\xi)
+\pi_{s}(\nap_X\xi)+\pi_{s+1}(\nap_X \xi))\\
\!\!\!&=&\!\!\!\phi_\theta(\pi_{s-1}(\nap_X\xi)+\pi_s(\nap_X\xi))
+S\circ\phi_\theta\circ\pi_{s+1}(\nabla_X^{\theta,\ell}\xi)\\
 \!\!\!&=&\!\!\!\psi_\theta(\nap_X\xi),
\eea
and this proves (\ref{pp}).

\medskip 

\noindent\emph{Case 2.}  For $ s= \ell $  we see that (\ref{s}) 
still holds. Using (\ref{one0}), (\ref{two0}), (\ref{D0}) and (\ref{cone}) 
we have
\bea
S(\nap_X\xi_\theta)
\!\!\!&=&\!\!\!S\circ\phi_\theta(\nabla_X^{\theta,\ell}\xi)
=S\circ\phi_\theta (\pi_{\ell-1}(\nabla_X^{\theta, \ell}\xi)
+\pi_{\ell}(\nabla_X^{\theta,\ell}\xi)+\pi_{\ell+1}(\nabla_X^{\theta,\ell}\xi))\\
 \!\!\!&=&\!\!\!\psi_\theta(\pi_{\ell-1}(\nap_X\xi)
+\pi_{\ell}(\nap_X\xi)+R^{\ell+1}_{-\theta}\pi_{\ell+1}(\nap_X R^{\ell}_\theta\xi))\\
 \!\!\!&=&\!\!\!\psi_\theta(\nap_X\xi),
\eea
and this proves (\ref{pp}).

\medskip 

\noindent\emph{Case 3.} Assume that $s=\ell+1$. 
Then, we have
\bea
\nap_XS\xi_\theta\!\!\!&=&\!\!\!\nap_X(\psi_\theta R^{\ell+1}_{-\theta} \xi)
=\psi_\theta(\nabla_X ^{\theta,\ell+r} R^{\ell+1}_{-\theta}\xi)\\
\!\!\!&=&\!\!\!\psi_\theta (\pi_{\ell}(\nap_X R^{\ell+1}_{-\theta}\xi)
+\pi_{\ell+1}(\nap_XR^{\ell+1}_{-\theta} \xi)+
\pi_{\ell+2}(\nap_X R^{\ell+1}_{-\theta}\xi)).
\eea
Using (\ref{one0}), (\ref{two0}), (\ref{D0}) and (\ref{cone}) 
we obtain
\bea
S(\nap_X\xi_\theta)\!\!\!&=&\!\!\! S\circ\phi_\theta(\nabla_X^{\theta,\ell}\xi)
=S\circ\phi_\theta(\pi_{\ell}(\nabla_X^{\theta,\ell}\xi)
+\pi_{\ell+1}(\nabla_X^{\theta,\ell}\xi)
+\pi_{\ell+2}(\nabla_X^{\theta,\ell}\xi))\\
\!\!\!&=&\!\!\!\psi_\theta(R^{\ell}_{-\theta}\pi_{\ell}(\nap_X\xi)
+R^{\ell+1}_{-\theta}\pi_{\ell+1}(\nap_X\xi)
+R^{\ell+2}_{-\theta}\pi_{\ell+2}(\nap_X\xi))\\
\!\!\!&=&\!\!\!\psi_\theta(\pi_{\ell}(\nap_X  R^{\ell+1}_{-\theta}\xi)
+\pi_{\ell+1}(\nap_XR^{\ell+1}_{-\theta}\xi)
+\pi_{\ell+2}(\nap_XR^{\ell+2}_{-\theta}\xi)),
\eea
and this proves (\ref{pp}).

\medskip 

\noindent\emph{Case 4.} Assume that $\ell+2\leq s\leq\ell+r+1$. 
Then, we have
\bea
\nap_XS\xi_\theta  
\!\!\!&=&\!\!\!\nap_X(\psi_\theta R^s_{-\theta} \xi)
=\psi_\theta(\nabla_X ^{\theta,\ell+r} R^s_{-\theta}\xi)\\
\!\!\!&=&\!\!\!\psi_\theta(\pi_{s-1}(\nap_XR^s_{-\theta}\xi)
+\pi_s(\nap_X R^s_{-\theta} \xi)
+\pi_{s+1}(\nap_X R^s_{-\theta}\xi)).
\eea
Using (\ref{one0}), (\ref{two0}), (\ref{D0}) and (\ref{cone}) 
we deduce that
\bea
S(\nap_X\xi_\theta)\!\!\!&=&\!\!\!S\circ\phi_\theta(\nabla_X^{\theta,\ell}\xi)
=S\circ\phi_\theta(\pi_{s-1}(\nabla_X^{\theta,\ell}\xi)
+\pi_s(\nabla_X^{\theta,\ell}\xi)+\pi_{s+1}(\nabla_X^{\theta,\ell}\xi))\\
\!\!\!&=&\!\!\!\psi_\theta (R^{s-1}_{-\theta}\pi_{s-1}(\nap_X\xi)
+R^s_{-\theta}\pi_s(\nap_X\xi)
+R^{s+1}_{-\theta}\pi_{s+1}(\nap_X\xi))\\
\!\!\!&=&\!\!\!\psi_\theta(\pi_{s-1}(\nap_XR^s_{-\theta}\xi)
+\pi_s( \nap_X  R^s_{-\theta}\xi)
+\pi_{s+1}(\nap_XR^s_{-\theta}\xi)),
\eea
and this proves (\ref{pp}).

\medskip 

\noindent\emph{Case 5.} Assume that $s=\ell+r$. Then, we have
\bea
\nap_XS\xi_\theta  
\!\!\!&=&\!\!\!\nap_X(\psi_\theta R^{\ell+r}_{-\theta}\xi)
=\psi_\theta (\nabla_X^{\theta,\ell+r}R^{\ell+r}_{-\theta}\xi)\\
\!\!\!&=&\!\!\!\psi_\theta(R^{\ell+r-1}_{-\theta}\pi_{\ell+r-1}(\nap_X\xi)
+\pi_{\ell+r}(\nap_XR^{\ell+r}_{-\theta}\xi)
+\pi_{\ell+r+1}(\nap_X R^{\ell+r}_{\theta}R^{\ell+r}_{-\theta}\xi)).
\eea
Using  (\ref{one0}), (\ref{two0}), (\ref{D0}) and (\ref{cone}) 
we obtain
\bea
S(\nap_X\xi_\theta)
\!\!\!&=&\!\!\!S\circ\phi_\theta (\nabla_X^{\theta,\ell}\xi)
=S\circ\phi_\theta(\pi_{\ell+r-1}(\nabla_X^{\theta,\ell}\xi)
+\pi_{\ell+r}(\nabla_X^{\theta,\ell}\xi)
+\pi_{\ell+r+1}(\nabla_X^{\theta,\ell}\xi))\\
\!\!\!&=&\!\!\!\psi_\theta(R^{\ell+r-1}_{-\theta}\pi_{\ell+r-1}(\nap_X\xi)
+R^{\ell+r}_{-\theta}\pi_{\ell+r}(\nap_X\xi)
+\pi_{\ell+r+1}(\nap_X\xi)),
\eea
and this proves (\ref{pp}).

\medskip 

\noindent\emph{Case 6.} Assume that $s=\ell+r+1$. Then, we have
\bea
\nap_XS\xi_\theta
\!\!\!&=&\!\!\!\nap_X\psi_\theta\xi
=\psi_\theta(\nabla_X^{\theta,\ell+r}\xi)\\
\!\!\!&=&\!\!\!\psi_\theta(R^{\ell+r}_{-\theta}\pi_{\ell+r}(\nap_X\xi)
+\pi_{\ell+r+1}(\nap_X\xi)+\pi_{\ell+r+2}(\nap_X\xi)).
\eea
Using (\ref{one0}), (\ref{two0}), (\ref{D0}) and (\ref{cone}) 
we obtain
\bea
S(\nap_X\xi_\theta)
\!\!\!&=&\!\!\!S\circ\phi_\theta(\nabla_X^{\theta,\ell}\xi)
=S\circ\phi_\theta(\pi_{\ell+r}(\nabla_X^{\theta,\ell}\xi)
+\pi_{\ell+r+1}(\nabla_X^{\theta,\ell}\xi)
+\pi_{\ell+r+2}(\nabla_X^{\theta,\ell}\xi))\\
\!\!\!&=&\!\!\!\psi_\theta(R^{\ell+r}_{-\theta}\pi_{\ell+r}(\nap_X\xi)
+\pi_{\ell+r+1}(\nap_X\xi)+\pi_{\ell+r+2}(\nap_X\xi)),
\eea
and this proves (\ref{pp}).

\medskip 

\noindent\emph{Case 7.} Assume that $s\geq\ell+r+2$. Then, we have
\bea
\nap_XS\xi_\theta  
\!\!\!&=&\!\!\!\nap_X\psi_\theta\xi
=\psi_\theta(\nabla_X ^{\theta,\ell+r}\xi)\\
\!\!\!&=&\!\!\!\psi_\theta(\pi_{s-1}(\nap_X\xi)+\pi_s(\nap_X\xi)
+\pi_{s+1}(\nap_X\xi)).
\eea
Using  (\ref{one0}), (\ref{two0}), (\ref{D0}) and (\ref{cone}) 
we obtain
\bea
S(\nap_X\xi_\theta)
\!\!\!&=&\!\!\!S\circ\phi_\theta(\nabla_X^{\theta,\ell}\xi)
=S\circ\phi_\theta(\pi_{s-1}(\nabla_X^{\theta,\ell}\xi)
+\pi_s(\nabla_X^{\theta,\ell}\xi)+\pi_{s+1}(\nabla_X^{\theta,\ell}\xi))\\
\!\!\!&=&\!\!\!\psi_\theta(\pi_{s-1}(\nap_X\xi)+\pi_s(\nap_X\xi)
+\pi_{s+1}(\nap_X\xi)),
\eea
and this proves (\ref{pp}).
Thus  $S$ is parallel, and the result follows. \qed

\section{Minimal submanifolds of rank 2}

Let $f\colon M^n\to\R^N$, $n\geq 3$, be a submanifold of rank two. This means
that the relative nullity subspaces $\Delta(x)\subset T_x M$
defined by
$$
\Delta (x) =\{X\in T_x M:\alpha_f(X,Y)=0\;\;\mbox{for all}\;\,Y\in T_x M\}
$$
form a codimension two subbundle of the tangent bundle. The submanifold is
called \emph{elliptic} if there exists an almost complex structure
$J\colon\Delta^\perp\to\Delta^\perp$ such that
$$
\alpha_f(X,X)+\alpha_f(JX,JX)=0\;\;\mbox{for all}\;\, X \in \Delta^{\perp}.
$$
Hence $f$ is minimal if and only if $J$ is orthogonal. As in the case of
elliptic surfaces, the normal bundle splits as
$$
N_fM=N^f_1\oplus \cdots \oplus N^f_{\tau}.
$$
It was shown in \cite{DF2}
that everything explained in this paper about polar surfaces to elliptic surfaces
extends to this case. In particular, any elliptic submanifold in case $(ii)$ admits
locally many polar surfaces which turn out to be elliptic surfaces.

Hereafter, we assume that  $f\colon M^n\to\R^N$  is minimal and simply connected
of rank two. For any $\varphi \in \Sp^1$ consider the tensor field 
$R_{\varphi}$ that is the identity on $\Delta$ and the  rotation through $\varphi$ 
in $\Delta^{\perp}$.
It was observed in \cite{DG} that the normal valued tensor field  given by
$$
\alpha_{\varphi}(X,Y)= \alpha_f(R_{\varphi}X,Y),
$$
satisfies the Gauss, Codazzi and Ricci equations with
respect to the normal connection of $f$. Hence, for each $\varphi\in\Sp^1$ there
exists a minimal submanifold $f_{\varphi}\colon M^n\to\R^N$ of rank two that forms 
the associated family of $f$.

According to the polar parametrization given in \cite[Thm.\ 10]{DF2}, minimal
submanifolds of rank two can be described parametrically along a subbundle of 
the normal bundle of an elliptic surface whose curvature ellipse of a specific 
order is circular. More precisely, given an elliptic surface 
$g\colon L^2\to\Q_c^{N-c}$, $c=0,1$, with $\mathcal{E}_\ell^g$ for some 
$1\leq \ell\leq\tau^o-1$ a circle, consider the  map 
$f\colon\Lambda_\ell\to\R^N$ defined by
$$
f(\delta)=h(x)+\delta, \;\; \delta \in \Lambda_\ell(x),
$$
where $\Lambda_\ell=N^g_{\ell+1}\oplus\cdots\oplus N^g_{\tau}$ and $h$ is any
$\ell$-cross
section to $g$, is at regular points a minimal submanifold of rank two with 
polar surface $g$. Conversely, any minimal submanifold of rank two admits 
locally such a parametrization with $g$ a polar map.

The recursive procedure for the construction of the cross sections
\cite[Prop.\ 6]{DF2} yields
$$
h=c\omega g + \grad\, \omega +\gamma_0+\gamma_1+\cdots+\gamma_\ell,
$$
where $\omega$ is a solution of the linear elliptic differential equation
$$
\Delta u+\<X,\grad u\> +c\lambda u =0
$$
for suitable $X \in TL$ and $\lambda\in C^\infty(L)$,  $\gamma_0$ is any section
in $\Lambda_\ell, \gamma_1 \in N^g_1$
is the unique solution of $A_{\gamma_1}= \mathrm{Hess}_{ \omega} +c \omega I$
and $\gamma_j \in N^g_j, 2 \leq j \leq \ell$,
where $L^2$ is endowed with the metric which makes $J$ orthogonal.

Take $g_{\theta}\in G_\ell$ and the corresponding vector 
bundle isometry $\phi_\theta\colon N_gL\to N_{g_\theta}L$.
Then,
$$
h_{\theta}=c\omega g_{\theta} + \grad\,\omega
+\phi_\theta\gamma_0+\phi_\theta\gamma_1+\cdots+\phi_\theta\gamma_\ell
$$
is an $\ell$-cross section to $g_{\theta} $. With these elements we have
the following result.

\begin{theorem}\po A submanifold in the associated family of a minimal
submanifold of rank two $f\colon M^n\to\R^N$  can be locally parametrized as
$$
f_{\theta}(\delta)=h_{\theta}(x)+\phi_\theta\delta, \;\; \delta \in
\Lambda_\ell(x).
$$
\end{theorem}

\proof It is easy to check that $f_{\theta}$ is isometric to $f$, has the same
normal connection and its second fundamental form is given by
$$
\alpha_{f_{\theta}}(X,Y)= \alpha_f(R_{-\theta}X,Y).
$$
In particular, $f_{\theta}$ is minimal and thus belongs to the associated family
of $f$.\qed\vspace{1,5ex}

The above discussion allows us to give an answer to the question of which minimal
submanifolds of rank two have  trivial  associated family. In fact, this
is the case if and only if the associated family of its polar surfaces is
trivial,
which is equivalent to the fact that a (local) bipolar surface to $f$, i.e.,
any polar surface to its polar surface,  is $m$-isotropic.

{\renewcommand{\baselinestretch}{1}

\hspace*{-20ex}\begin{tabbing} \indent\= IMPA -- Estrada Dona Castorina, 110
\indent\indent\= Univ. of Ioannina -- Math. Dept. \\
\> 22460-320 -- Rio de Janeiro -- Brazil  \>
45110 Ioannina -- Greece \\
\> E-mail: marcos@impa.br \> E-mail: tvlachos@uoi.gr
\end{tabbing}}

\end{document}